\documentclass[11pt]{article}
\usepackage{amsmath}
\usepackage{pxfonts}
\usepackage{yfonts}
\usepackage{dsfont}
\usepackage{marvosym}
\usepackage{graphicx}
\usepackage{relsize}

\def\bel{\begin{equation}\label}
\def\eeq{\end{equation}}
\def\ds{\displaystyle}
\def\endproof{\hphantom{MM}
\hfill\llap{$\square$}\goodbreak}

\def\mt{\longrightarrow}
\def\v{\vskip 1em}

\def\R{\mathds R}

\def\C{\mathfrak{C}}
\def\Cx{\mathds C}
\def\N{{\bf N}}

\def\Re{{\bf Re}}

\def\S{{\bf S}}

\def\J{{\bf J}}

\def\L{{\bf L}}

\def\i{{\bf i}}
\def\Tilde{\widetilde}
\def\Hat{\widehat}

\def\bar{\overline}

\def\I{{\bf I}}

\def\alpha{\alphaup}
\def\beta{\betaup}
\def\gamma{\gammaup}
\def\delta{\deltaup}
\def\xi{{\xiup}}
\def\eta{{\etaup}}
\def\tau{{\tauup}}
\def\rho{{\rhoup}}
\def\phi{{\phiup}}
\def\psi{{\psiup}}
\def\lambda{{\lambdaup}}
\def\omega{\omegaup}
\def\varphi{{\varphiup}}
\def\gamma{{\gammaup}}

\def\a{{\bf a}}

\def\e{{\bf e}}

\begin{document}
 \[\begin{array}{cc}\hbox{\LARGE{\bf Fractional integration with singularity on Light-cone, II}}
  \end{array}\]
 
  \[\hbox{Zipeng Wang}\]
 
 \begin{abstract}
 We study a family of convolution operators whose kernels carrying a critical index have singularity appeared on the light-cone in $\R^{n+1}$. We improve an result obtained by 
 Oberlin \cite{Oberlin} in 1989.  \end{abstract}
 
 \section{Introduction}
 \setcounter{equation}{0}
Let $\alpha,v\in\R$ and
\bel{lambda-alpha}
\lambda(\alpha)~=~{n+1\over 2}\left(1-{\alpha\over n}\right)+\i v,\qquad \gamma(\alpha,v)~=~\pi^{-\lambda(\alpha)}\Gamma^{-1}\left(1-\lambda(\alpha)\right)
\eeq
where $\Gamma$ is Gamma function.

For $0<\Re\lambda(\alpha)<1$, we define
\bel{Omega^alpha}
\begin{array}{ccc}\ds
\Omega^{\lambda(\alpha)}(x)~=~\left\{\begin{array}{lr}\ds
\gamma(\alpha,v)  \left({1\over 1-|x|^2}\right)^{\lambda(\alpha)}, \qquad |x|<1,
\\\\ \ds~~~~~~~~~~~~~
0,\qquad\qquad~~~~~~~~~~~|x|\ge1
\end{array}\right.
\end{array}
\eeq
  whose Fourier transform equals
 \bel{Omega^alpha Transform} 
\begin{array}{cc}\ds
\Hat{\Omega}^{\lambda(\alpha)}(\xi)~=~\left({1\over|\xi|}\right)^{{n\over 2}-\lambda(\alpha)} \J_{{n\over 2}-\lambda(\alpha)}\Big(2\pi|\xi|\Big)
\end{array}
\eeq
where $\J$ is Bessel function. See chapter IV of Stein and Weiss \cite{Stein-Weiss*}.

Denote  $\gamma(\alpha)=\gamma(\alpha,0)$ and $\Omega^\alpha=\Omega^{\lambda(\alpha)}$ at $v=0$. 
Observe that  
\bel{Cone singular}
 \begin{array}{lr}\ds
~\gamma(\alpha) \iint_{|u|<|r|} f(x-u,t-r)  \left({1\over r^2-|u|^2}\right)^{{n+1\over 2}\left(1-{\alpha\over n}\right)}dudr
\\\\ \ds~
~=~\iint_{\R^{n+1}} f(x-u,t-r)\Omega^\alpha\left({u\over r}\right)|r|^{\left({n+1\over n}\right)\alpha-1-n}dudr.
\end{array}
 \eeq 
 $\Omega^{\lambda(\alpha)}$ can be extensively defined for $0<\alpha< n$ by the inverse Fourier transform of (\ref{Omega^alpha Transform}), or equivalently by an analytic continuation from (\ref{Omega^alpha}).\footnote{See p.462 in  \cite{Strichartz}.}
 
 Let $\Omega^{\lambda(\alpha)}_r, r>0$   be a dilate of $\Omega^{\lambda(\alpha)}$ in the sense of distribution.  
A fractional integral operator $\I_\alpha$ is defined by
 \bel{I_alpha}
 \begin{array}{ccc}\ds
\I_\alpha f(x,t)~=~ \iint_{\R^{n+1}} f(x-u,t-r)\Omega^\alpha_{|r|}(u)| r|^{\left({n+1\over n}\right)\alpha-1}dudr,\qquad 0<\alpha<n.
\end{array}
\eeq 
$\diamond$ We always write $\C$ as a generic constant depending on its subindices.

 {\bf Conjecture A} ~~{\it Let $\I_\alpha$ defined in (\ref{I_alpha}) for $0<\alpha<n$. We have
\bel{RESULT}
 \left\| \I_\alpha f\right\|_{\L^q(\R^{n+1})}~\leq~\C_{p~q}~\left\|f\right\|_{\L^p(\R^{n+1})},\qquad 1<p<q<\infty
\eeq
 if and only if
\bel{FORMULA+CONSTRAINT}
\begin{array}{cc}\ds
 {\alpha\over n}~=~{1\over p}-{1\over q},
\qquad
 {n-1\over 2n}+\left({n+1\over 2n}\right){\alpha\over n}~<~{1\over p}~<~{n+1\over 2n}+\left({n-1\over 2n}\right){\alpha\over n}.
 \end{array}
 \eeq}

The necessity of constraints (\ref{FORMULA+CONSTRAINT}) has been known for the norm inequality  (\ref{RESULT}) to hold.   The converse
 is proved for ${n\over n+1}\leq \alpha<n$ by Oberlin \cite{Oberlin} who proposed the question of the analogue result when $0<\alpha<{n\over n+1}$.  See the region {\bf I} in Figure 1 below. In this paper,  we prove the desired regularity of $\I_\alpha$ for $\left({1\over p},\alpha\right)$ inside the region {\bf II} .
 
 \begin{figure}[h]
\centering
\includegraphics[scale=0.52]{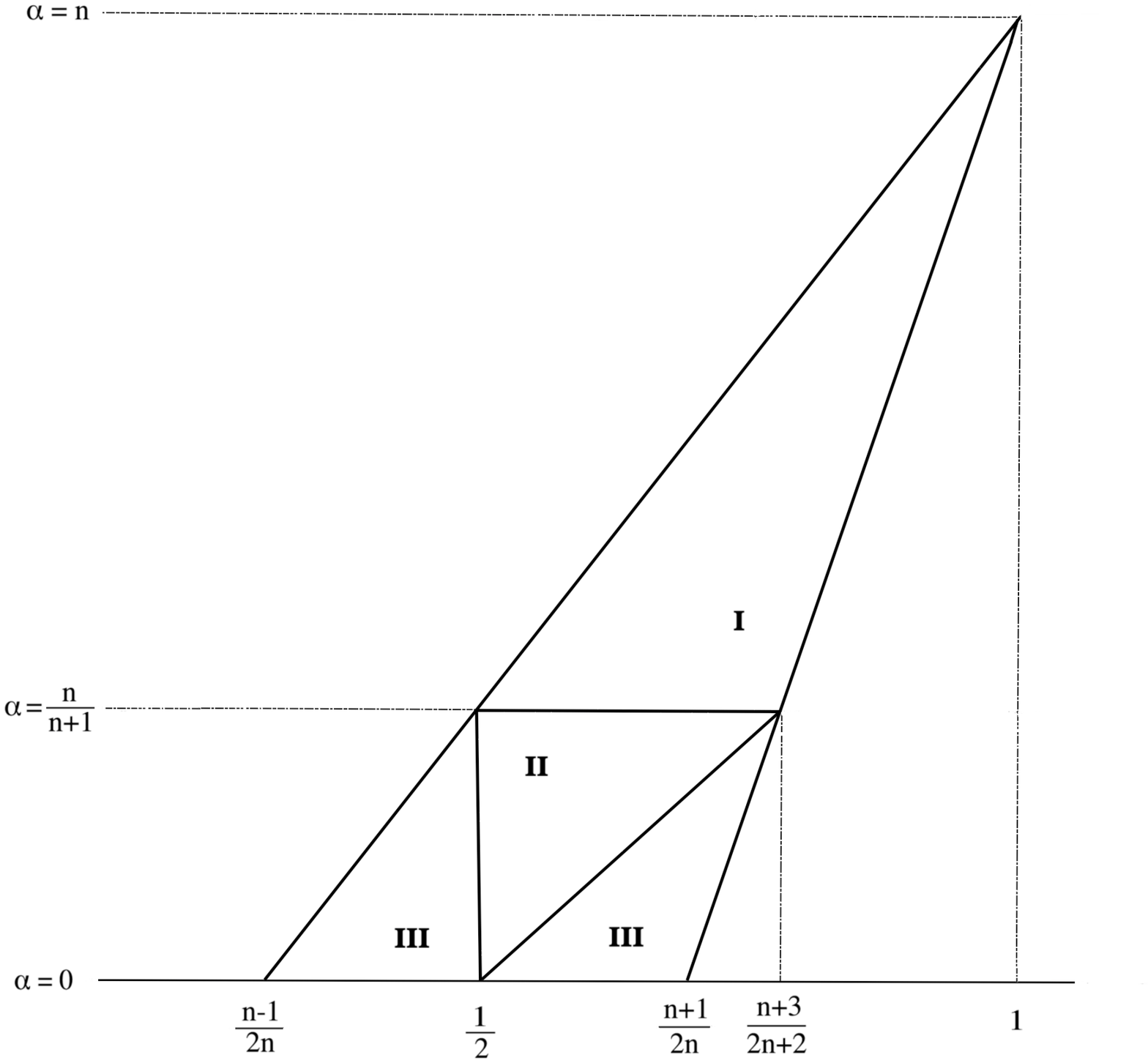}
\caption{The horizontal direction measures the range of ${1\over p}$.}
\end{figure} 

Some other regarding works  refer to Ricci and Stein \cite{Ricci-Stein}, Greenleaf \cite{Greenleaf} and Littman \cite{Littman}. 
Historical background  can be found in the book by Gelfand and Shilov \cite{Gelfand-Shilov}.

\section{Statement of main results}
\setcounter{equation}{0}

{\bf Theorem One} ~~{\it Let $\I_\alpha$ defined in (\ref{I_alpha}) for $0<\alpha\leq{n\over n+1}$. We have
\bel{RESULT*}
\begin{array}{cc}\ds
 \left\| \I_\alpha f\right\|_{\L^q(\R^{n+1})}~\leq~\C_{p~q}~\left\|f\right\|_{\L^p(\R^{n+1})},\qquad 1<p<q<\infty
\\\\ \ds
\hbox{if}\qquad {\alpha\over n}~=~{1\over p}-{1\over q},\qquad
{1\over 2}~<~{1\over p}~<~{1\over 2}+{\alpha\over n}.
 \end{array}
 \eeq}

 {\bf Sketch of proof:~} We follow the lines from the elegant work by Oberlin \cite{Oberlin} but  improve
 an important lemma that plays a principal role in our analysis. This is done in section 3 where the remaining section is devoted to certain regularity estimates $w.~r.~t$ the line of duality.   Consequently, we obtain the following result.
 
{\bf Proposition 1}~~ {\it Let $0<\alpha\leq{n\over n+1}$ and $\Omega^\alpha$ defined by the inverse Fourier transform of  (\ref{Omega^alpha Transform}). We have
\bel{Result one}
\begin{array}{cc}\ds
\left\{ \int_0^\infty r^{\alpha q} \left\{ \int_{\R^n} \left|\Big(f\ast\Omega^{\alpha}_r\Big)(x)\right|^qdx\right\}{dr\over r}\right\}^{1\over q} ~\leq~\C_{\alpha}~\left\|f\right\|_{\L^{2}(\R^n)}
\\\\ \ds
\hbox{for}\qquad {\alpha\over n}~=~{1\over 2}-{1\over q},\qquad 1<q<\infty.
 \end{array}
 \eeq}
 
 On the other hand, we borrow the next result from Oberlin \cite{Oberlin}  ( See (19) on p. 754 ).
  
{\bf Proposition 2}~~ {\it Let $0<\alpha<n$ and $\Omega^{\lambda(\alpha)}$ defined by the inverse Fourier transform of  (\ref{Omega^alpha Transform}).  We have
\bel{convolution inequality}
\begin{array}{cc}\ds
 \left\|  f\ast\Omega^{\lambda(\alpha)}\right\|_{\L^q(\R^n)}~\leq~\C_{\alpha~v}~\left\|f\right\|_{\L^p(\R^n)}
 \\\\ \ds
\hbox{for}\qquad  {1\over p}~=~{1\over 2}\left(1+{\alpha\over n}\right),\qquad {1\over q}~=~{1\over 2}\left(1-{\alpha\over n}\right).
 \end{array}
\eeq}

In section 4, by using {\bf Proposition 1} and {\bf 2} in a combination of one dimensional fractional integration with convolution inequalities on $\R^n$, \footnote{the original idea of this argument is given by  Strichartz \cite{Strichartz'}-\cite{Strichartz}.} we complete the proof of {\bf Theorem One}.

 \section{Regularity estimate $w.r.t$ the line of duality}
 \setcounter{equation}{0}
  Let  $\Omega^{z}, z\in\Cx$ be a distribution on $\R^n$ whose dilate $\Omega^z_r,~r>0$ and adjoint $\Tilde{\Omega}^z$ are defined as
  \bel{dila adjoint}
  \int_{\R^n} \phi(x) d\Omega^z_r(x)~=~\int_{\R^n} \phi(rx)d\Omega^z(x),\qquad \int_{\R^n}\phi(x)d\Tilde{\Omega}^z(x)~=~\int_{\R^n}\phi(-x)d\bar{\Omega^z}(x)
  \eeq
  for all test functions $\phi$.

 {\bf Lemma One} ~~~~{\it
  Let $0<\alpha\leq{n\over n+1}$.  We have
\bel{Key}
\begin{array}{cc}\ds
\left\{ \int_0^\infty r^{\alpha q} \left\{ \int_{\R^n} \left|\Big(f\ast\Omega^{z}_r\Big)(x)\right|^qdx\right\}{dr\over r}\right\}^{1\over q} ~\leq~\C_{\alpha~z}~\left\|f\right\|_{\L^{2}(\R^n)}
\\\\ \ds
 \hbox{for}\qquad {\alpha\over n}~=~{1\over 2}-{1\over q},\qquad 1<q<\infty
 \end{array}
  \eeq
  if
   \bel{crucial est}
 \left\| f\ast\Omega^{z}_{r}\ast\Tilde{\Omega}^{z}_s \right\|_{\L^q(\R^{n})}~\leq~\C_{\alpha~z}~\left({1\over r}\right)^{\left({n+1\over 2n}\right)\alpha}\left({1\over |r-s|}\right)^{\left({n-1\over n}\right)\alpha} \left({1\over s}\right)^{\left({n+1\over 2n}\right)\alpha}\left\| f\right\|_{\L^{q\over q-1}(\R^{n})}
 \eeq
 for every $r,s>0$. }
\v

{\bf Proof:} Consider
\bel{Sf} 
 \Big(\S f\Big)(x,r) ~=~|r|^{\alpha-{1\over q}} \Big(f\ast\Omega^{z}_{|r|}\Big)(x).
 \eeq
 Observe that (\ref{Key}) is equivalent to
 \bel{Result 2,q}
 \left\| \S f\right\|_{\L^q(\R^{n+1})}~\leq~\C_{\alpha~z}~\left\|f\right\|_{\L^2(\R^n)}. 
 \eeq
Let $\S^*$ denote the adjoint operator of $\S$. We  prove (\ref{Result 2,q}) by showing $\S^*\colon\L^{q\over q-1}(\R^{n+1})\mt\L^2(\R^n)$. This in turn can be deduced from 
$\S\S^*\colon\L^{q\over q-1}(\R^{n+1})\mt\L^q(\R^{n+1})$. 
 
Let $g\colon\R^{n+1}\mt\R$. From direct computation, we have
\bel{SS*g}
\Big(\S\S^* g\Big)(x,r)~=~ \int_\R |rs|^{\alpha-{1\over q}}\Big(g(\cdot,s)\ast\Omega^{z}_{|r|}\ast\Tilde{\Omega}^{z}_{|s|}\Big)(x) ds.
 \eeq
 Set
 \bel{gamma}
 \begin{array}{cc}\ds
 -\gamma~=~\alpha-{1\over q}-\left({n+1\over 2n}\right)\alpha~=~ \left({n-1\over 2n}\right)\alpha-{1\over q}.
 \end{array} 
\eeq
Note that ${\alpha\over n}={1\over 2}-{1\over q}$ and $0<\alpha\leq{n\over n+1}$ together imply
 \bel{power range}
 \begin{array}{lr}\ds
  \left({n-1\over 2n}\right)\alpha-{1\over q}~=~ \left({n-1\over 2n}\right)\alpha-{1\over 2}+{\alpha\over n}
  \\\\ \ds~~~~~~~~~~~~~~~~~~~~~~~
 ~=~\left({n+1\over 2n}\right)\alpha-{1\over 2} ~\leq~0.
 \end{array}
 \eeq 
 From (\ref{gamma})-(\ref{power range}), we have
 \bel{gamma range}
 \begin{array}{cc}\ds
 0~\leq~ \gamma~=~{1\over 2}-\left({n+1\over 2n}\right)\alpha~=~{1\over q}- \left({n-1\over 2n}\right)\alpha~<~{1\over q}.
 \end{array}
\eeq
Moreover, 
\bel{index computa}
\begin{array}{lr}\ds
1-2\left({\alpha\over n}+\gamma\right)~=~1-2\left[{\alpha\over n}+{1\over 2}-\left({n+1\over 2n}\right)\alpha\right]~=~\left({n-1\over n}\right)\alpha.
\end{array}
\eeq

Recall a classical result obtained by Stein and Weiss \cite{Stein-Weiss}.

{\bf Stein-Weiss Theorem (1958)}  ~{\it Let $0<\a<\N$ and $\gamma, \delta\in\R$. We have
\bel{norm ineq}
\left\{\int_{\R^\N}\left|\int_{\R^\N} f(u)\left({1\over|u|}\right)^\delta\left({1\over |x-u|}\right)^{\N-\a}\left({1\over |x|}\right)^{\gamma } du\right|^q  dx\right\}^{1\over q}~\leq~\C_{p~q~\gamma~\delta~a} \left\{\int_{\R^\N} \left|f(x)\right|^p  dx\right\}^{1\over p}
\eeq
for  $1<p\leq q<\infty$ if 
\bel{constraints powers}
\begin{array}{cc}\ds
\gamma~<~{\N\over q},\qquad \delta~<~\N\left({p-1\over p}\right),\qquad \gamma+\delta~\ge~0
\end{array}
\eeq
and
\bel{Formula*}
{\a\over \N}~=~{1\over p}-{1\over q}+{\gamma+\delta\over \N}.
\eeq}

From (\ref{SS*g}), we have
 \bel{SS*g Est} 
 \begin{array}{lr}\ds
 \left\| \S\S^* g\right\|_{\L^q(\R^{n+1})} 
 ~=~\left\{ \iint_{\R^{n+1}}\left|\int_\R |rs|^{\alpha-{1\over q}}\Big(g(\cdot,s)\ast\Omega^{z}_{|r|}\ast\Tilde{\Omega}^{z}_{|s|}\Big)(x) ds \right|^q dxdr\right\}^{1\over q}
 \\\\ \ds~~~~~~~
 ~\leq~\left\{\int_\R\left\{\int_\R |rs|^{\alpha-{1\over q}}\left\{\int_{\R^{n}}\left|\Big(g(\cdot,s)\ast\Omega^{z}_{|r|}\ast\Tilde{\Omega}^{z}_{|s|}\Big)(x)\right|^q dx\right\}^{1\over q} ds \right\}^q dr\right\}^{1\over q}
~~
 \hbox{\small{by Minkowski integral inequality}}
 \\\\ \ds ~~~~~~~
  ~\leq~\C_{\alpha~z} \left\{\int_\R\left\{\int_\R |rs|^{\alpha-{1\over q}}\left({1\over |rs|}\right)^{\left({n+1\over 2n}\right)\alpha}\left({1\over |r-s|}\right)^{\left({n-1\over n}\right)\alpha} \left\| g(\cdot,s)\right\|_{\L^{q\over q-1}(\R^{n})} ds \right\}^q dr\right\}^{1\over q}
\qquad
 \hbox{\small{by (\ref{crucial est})}}  
  \\\\ \ds ~~~~~~~
 ~=~\C_{\alpha~z} \left\{\int_\R \left\{\int_\R \left({1\over |r|}\right)^{\gamma}\left({1\over |r-s|}\right)^{1-2\left({\alpha\over n}+\gamma\right)}\left({1\over|s|}\right)^\gamma \left\| g(\cdot,s)\right\|_{\L^{q\over q-1}(\R^{n})} ds \right\}^q dr\right\}^{1\over q} 
~~~~\hbox{\small{by (\ref{gamma})-(\ref{index computa})}}
\\\\ \ds ~~~~~~~
~\leq~\C_{\alpha~z} \left\{\int_\R \left\| g(\cdot,r)\right\|_{\L^{q\over q-1}(\R^n)}^{q\over q-1} dr\right\}^{q-1\over q}\qquad\hbox{\small{by {\bf Stein-Weiss theorem}}}
 \end{array}
 \eeq
 where the last inequality is obtained by taking into account that
$\N=1,~ p={q\over q-1},~ \delta=\gamma$  and $\a=2\left({\alpha\over n}-\gamma\right)$ in (\ref{constraints powers})-(\ref{Formula*}).\endproof
 \v
Recall some well known estimates of Bessel functions: 
For $\mu>-{1\over 2}$ and $\rho>0$, 
\bel{Bessel formula}
 \J_{\mu+\i v}(\rho)~=~\left({2\over \pi \rho}\right)^{1\over 2}\cos\left(\rho-{\pi\over 2}(\mu+\i v)-{\pi\over 4}\right)+\e_v(\rho),
\eeq
\bel{Bessel norm}
 |\e_v(\rho)|~\leq~\C_\mu~ e^{\pi|v|}\left\{\begin{array}{lr}\ds \rho^{-{1\over 2}},\qquad 0<\rho\leq1,
 \\ \ds
 \rho^{-{3\over 2}},\qquad~~~ \rho>1.
 \end{array}\right.
\eeq
 In particular,  (\ref{Bessel formula})-(\ref{Bessel norm}) can be obtained by carrying out  a complex analogue to the proof of Lemma 3.11 on {p.158}  in Stein and Weiss \cite{Stein-Weiss*}.  More discussion can be found in  Watson \cite{Watson}.

Let $\Omega^{\lambda(\alpha)}$  defined by the inverse Fourier transform of (\ref{Omega^alpha Transform}). Consider
 \bel{Omega split}
 \Hat{\Omega}^{\lambda(\alpha)}(\xi)~=~\Hat{\mathfrak{S}}^{\lambda(\alpha)}(\xi)+\Hat{\mathcal{E}}^{\lambda(\alpha)}(\xi)
 \eeq
 where
 \bel{Omega split S+E}
 \begin{array}{cc}\ds
\Hat{\mathfrak{S}}^{\lambda(\alpha)}(\xi)~=~{1\over \pi}  \left({1\over |\xi|}\right)^{\left({n+1\over 2n}\right)\alpha+\i v} \cos\left(2\pi|\xi|-{\pi\over 2}\left({n+1\over 2n}\right)\alpha-{\pi\over 2}(1+\i v)\right), 
 \\\\ \ds
 \Hat{\mathcal{E}}^{\lambda(\alpha)}(\xi)~=~  \left({1\over |\xi|}\right)^{\left({n+1\over 2n}\right)\alpha-{1\over 2}+\i v}\e_v(2\pi|\xi|). 
 \end{array}
 \eeq 
Let $\Omega^\alpha=\Omega^{\lambda(\alpha)},~v=0$. We aim to prove (\ref{Key}) for $f\ast\Omega^\alpha_r,~r>0$ by showing 
 \bel{crucial est.1}
 \left\| f\ast\mathfrak{S}^{\alpha}_r\ast\Tilde{\mathfrak{S}}^{\alpha}_s \right\|_{\L^q(\R^{n})}~\leq~\C_{\alpha}~\left({1\over rs}\right)^{\left({n+1\over 2n}\right)\alpha}\left({1\over |r-s|}\right)^{\left({n-1\over n}\right)\alpha} \left\| f\right\|_{\L^{q\over q-1}(\R^{n})}
 \eeq
 and
  \bel{crucial est.2}
 \left\| f\ast\mathcal{E}^{\alpha}_r\ast\Tilde{\mathcal{E}}^{\alpha}_s\right\|_{\L^q(\R^{n})}~\leq~\C_{\alpha}~\left({1\over rs}\right)^{\left({n+1\over 2n}\right)\alpha}\left({1\over |r-s|}\right)^{\left({n-1\over n}\right)\alpha} \left\| f\right\|_{\L^{q\over q-1}(\R^{n})}
 \eeq
 for every $r,s>0$. 
 
 From (\ref{Omega split S+E}), we have
 \bel{f E E}
 \begin{array}{lr}\ds
\Big( f\ast\mathcal{E}^{\alpha}_r\ast\Tilde{\mathcal{E}}^{\alpha}_s \Big)(x)
~=~ \int_{\R^n} e^{2\pi\i x\cdot\xi} \Hat{f}(\xi) \left({1\over r|\xi|}\right)^{\left({n+1\over 2n}\right)\alpha-{1\over 2}}\e_o(2\pi r|\xi|)  \left({1\over s|\xi|}\right)^{\left({n+1\over 2n}\right)\alpha-{1\over 2}}\bar{\e_o}(2\pi s|\xi|) d\xi
\\\\ \ds~~~~~~~~~~~~~~~~~
~=~  \left({1\over rs}\right)^{\left({n+1\over 2n}\right)\alpha}\int_{\R^n} e^{2\pi\i x\cdot\xi} \Hat{f}(\xi)\left({1\over |\xi|}\right)^{\left({n+1\over n}\right)\alpha}\left(r|\xi|\right)^{{1\over 2}}\left(s|\xi|\right)^{{1\over 2}}\e_o(2\pi r|\xi|) \bar{\e_o}(2\pi s|\xi|) d\xi.
\end{array}
\eeq 
Without lose of the generality,  assume $r\ge s$ so that $r^{-1}\leq(r-s)^{-1}$.

{\bf Case 1}~~ Suppose $|\xi|\leq (2\pi r)^{-1}$.
 We have
 \bel{E norm est 1}
 \begin{array}{lr}\ds
\left({1\over |\xi|}\right)^{\left({n+1\over n}\right)\alpha}\left(r|\xi|\right)^{1\over 2}\left(s|\xi|\right)^{1\over 2}\left|\e_o(2\pi r|\xi|) \bar{\e_o}(2\pi s|\xi|)\right|
\\\\ \ds
~\leq~\C \left({1\over |\xi|}\right)^{\left({n+1\over n}\right)\alpha}~\leq~
  \C \left({1\over r}\right)^{\left({n-1\over n}\right)\alpha}\left({1\over |\xi|}\right)^{2\alpha}\qquad\hbox{\small{by (\ref{Bessel formula})}}
  \\\\ \ds
  ~\leq~\C \left({1\over r-s}\right)^{\left({n-1\over n}\right)\alpha}\left({1\over |\xi|}\right)^{2\alpha}.  
 \end{array}
 \eeq
 
 {\bf Case 2}~~ Suppose $|\xi|> (2\pi r)^{-1}$ and $|\xi|\leq (2\pi s)^{-1}$. We have 
  \bel{E norm est 2}
 \begin{array}{lr}\ds
 \left({1\over |\xi|}\right)^{\left({n+1\over n}\right)\alpha}\left(r|\xi|\right)^{1\over 2}\left(s|\xi|\right)^{1\over 2}\left|\e_o(2\pi r|\xi|) \bar{\e_o}(2\pi s|\xi|)\right|
\\\\ \ds
  ~\leq~\C  \left({1\over r}\right)\left({1\over |\xi|}\right)^{\left({n+1\over n}\right)\alpha+1} ~\leq~\C \left({1\over r}\right)^{\left({n-1\over n}\right)\alpha}\left({1\over |\xi|}\right)^{2\alpha}     \qquad\hbox{\small{by (\ref{Bessel formula})}}
      \\\\ \ds
  ~\leq~ \C \left({1\over r-s}\right)^{\left({n-1\over n}\right)\alpha}\left({1\over |\xi|}\right)^{2\alpha}.
 \end{array}
 \eeq
 
{\bf Case 3}~~ Suppose $|\xi|> (2\pi s)^{-1}$. We have
 \bel{E norm est 3}
 \begin{array}{lr}\ds
 \left({1\over |\xi|}\right)^{\left({n+1\over n}\right)\alpha}\left(r|\xi|\right)^{1\over 2}\left(s|\xi|\right)^{1\over 2}\left|\e_o(2\pi r|\xi|) \bar{\e_o}(2\pi s|\xi|)\right|
\\\\ \ds
  ~\leq~\C  \left({1\over rs}\right)\left({1\over |\xi|}\right)^{\left({n+1\over n}\right)\alpha+2} ~\leq~\C \left({1\over r}\right)\left({1\over |\xi|}\right)^{\left({n+1\over n}\right)\alpha+1}  \qquad\hbox{\small{by (\ref{Bessel formula})}}
   \\\\ \ds
  ~\leq~\C \left({1\over r}\right)^{\left({n-1\over n}\right)\alpha}\left({1\over |\xi|}\right)^{2\alpha} 
 ~\leq~ \C \left({1\over r-s}\right)^{\left({n-1\over n}\right)\alpha}\left({1\over |\xi|}\right)^{2\alpha}.
 \end{array}
 \eeq
Because ${2\alpha\over n}={q-1\over q}-{1\over q}$,   note that $|\xi|^{-2\alpha}$ is a $\L^{q\over q-1}(\R^n)\mt\L^q(\R^n)$ Fourier multiplier. Indeed, its inverse Fourier transform equals $\C_\alpha |x|^{2\alpha-n}$.
  Furthermore, ${q-1\over q}-{1\over 2}={\alpha\over n}={1\over 2}-{1\over q}$. Hence that $|\xi|^{-\alpha}$ is simultaneously a $\L^{q\over q-1}(\R^n)\mt\L^2(\R^n)$ or $\L^2(\R^n)\mt\L^q(\R^n)$ Fourier multiplier.
  
From (\ref{E norm est 1})-(\ref{E norm est 3}),  Hardy-Littlewood-Sobolev inequality \cite{Hardy-Littlewood}-\cite{Sobolev} together with Plancherel theorem imply
  \bel{f E E result}
  \begin{array}{lr}\ds
 \left\{ \int_{\R^n} \left|\Big( f\ast\mathcal{E}^{\alpha}_r\ast\Tilde{\mathcal{E}}^{\alpha}_s \Big)(x)  \right|^q dx\right\}^{1\over q}
 \\\\ \ds
 ~\leq~\C_{\alpha} \left({1\over rs}\right)^{\left({n+1\over 2n}\right)\alpha}   \left({1\over r-s}\right)^{\left({n-1\over n}\right)\alpha}\left\{\int_{\R^n}\left|f(x)\right|^{q\over q-1}dx\right\}^{q-1\over q}.
 \end{array}
 \eeq
 Return to (\ref{crucial est.1}). From direct computation,  we have 
 \bel{f S S}
 \begin{array}{lr}\ds
\Big( f\ast\mathfrak{S}^{\alpha}_r\ast\Tilde{\mathfrak{S}}^{\alpha}_s \Big)(x)
~=~{1\over \pi^2} \int_{\R^n} e^{2\pi\i x\cdot\xi} \Hat{f}(\xi) \left({1\over r|\xi|}\right)^{\left({n+1\over 2n}\right)\alpha} \cos\left(2\pi r|\xi|-{\pi\over 2}\left({n+1\over 2n}\right)\alpha-{\pi\over 2}\right)
\\\\ \ds~~~~~~~~~~~~~~~~~~~~~~~~~~~~~~~~~~~~~
 \left({1\over s|\xi|}\right)^{\left({n+1\over 2n}\right)\alpha} \cos\left(2\pi s|\xi|-{\pi\over 2}\left({n+1\over 2n}\right)\alpha-{\pi\over 2}\right) d\xi
 \\\\ \ds~~~~~~~~~~~~~~~
 ~=~\C_{\alpha} \left({1\over rs}\right)^{\left({n+1\over 2n}\right)\alpha}\int_{\R^n} e^{2\pi\i x\cdot\xi} \Hat{f}(\xi) \left({1\over |\xi|}\right)^{\left({n+1\over n}\right)\alpha}\cos\left(2\pi(r+s)|\xi|+{\pi\over 2}\left({n+1\over n}\right)\alpha\right)d\xi
  \\\\ \ds~~~~~~~~~~~~~~~
 ~+~\C_{\alpha} \left({1\over rs}\right)^{\left({n+1\over 2n}\right)\alpha}\int_{\R^n} e^{2\pi\i x\cdot\xi} \Hat{f}(\xi) \left({1\over |\xi|}\right)^{\left({n+1\over n}\right)\alpha}\cos\left(2\pi(r-s)|\xi|+{\pi\over 2}\left({n+1\over n}\right)\alpha\right)d\xi
  \\\\ \ds~~~~~~~~~~~~~~~
 ~+~\C_{\alpha} \left({1\over rs}\right)^{\left({n+1\over 2n}\right)\alpha}\int_{\R^n} e^{2\pi\i x\cdot\xi} \Hat{f}(\xi) \left({1\over |\xi|}\right)^{\left({n+1\over n}\right)\alpha}\sin\left(2\pi(r-s)|\xi|+{\pi\over 2}\left({n+1\over n}\right)\alpha\right)d\xi. 
 \end{array}
 \eeq
 From (\ref{f S S}), it is suffice to show that
 \bel{Fourier multiplier r-s}
\left({1\over |\xi|}\right)^{\left({n+1\over n}\right)\alpha}\sin\left(2\pi|\xi|+{\pi\over 2}\left({n+1\over n}\right)\alpha\right),\qquad 
 \left({1\over |\xi|}\right)^{\left({n+1\over n}\right)\alpha}\cos\left(2\pi|\xi|+{\pi\over 2}\left({n+1\over n}\right)\alpha\right)
 \eeq 
 are $\L^{q\over q-1}(\R^n)\mt\L^q(\R^n)$ Fourier multipliers. 
 
 We then have
 \bel{f S cos Est}
 \begin{array}{lr}\ds
\left\{\int_{\R^n} \left|  \left({1\over rs}\right)^{\left({n+1\over 2n}\right)\alpha}\int_{\R^n} e^{2\pi\i x\cdot\xi} \Hat{f}(\xi) \left({1\over |\xi|}\right)^{\left({n+1\over n}\right)\alpha}\cos\left(2\pi(r-s)|\xi|+{\pi\over 2}\left({n+1\over n}\right)\alpha\right)d\xi\right|^q dx \right\}^{1\over q}  
 \\\\ \ds
~=~\left({1\over rs}\right)^{\left({n+1\over 2n}\right)\alpha}(r-s)^{\left({n+1\over n}\right)\alpha}  \left\{\int_{\R^n} \left| \int_{\R^n} e^{2\pi\i x\cdot\xi} \Hat{f}(\xi) \left({1\over (r-s)|\xi|}\right)^{\left({n+1\over n}\right)\alpha}\cos\left(2\pi(r-s)|\xi|+{\pi\over 2}\left({n+1\over n}\right)\alpha\right)d\xi\right|^q dx \right\}^{1\over q} 
\\\\ \ds
~=~\left({1\over rs}\right)^{\left({n+1\over 2n}\right)\alpha}(r-s)^{\left({n+1\over n}\right)\alpha} \left\{\int_{\R^n} \left| \int_{\R^n} e^{2\pi\i x\cdot\xi} \Hat{f}((r-s)^{-1}\xi) \left({1\over |\xi|}\right)^{\left({n+1\over n}\right)\alpha}\cos\left(2\pi|\xi|+{\pi\over 2}\left({n+1\over n}\right)\alpha\right)d\xi\right|^q dx \right\}^{1\over q}  
\\ \ds~~~~~~~~~~~~~~~~~~~~~~~~~~~~~~~~~~~~~~~~~~~~~~~~~~~~~~~~~~~~~~~~~~~~~~~~~~~\xi\mt(r-s)^{-1}\xi,\qquad x\mt (r-s) x
\\ \ds
~=~\left({1\over rs}\right)^{\left({n+1\over 2n}\right)\alpha}(r-s)^{\left({n+1\over n}\right)\alpha+{n\over q}} 
\\ \ds~~~~~~~~~~~~~~~~~~~~~~~~
\left\{\int_{\R^n} \left| \int_{\R^n} f((r-s)u)\left\{\int_{\R^n} e^{2\pi\i (x-u)\cdot\xi}  \left({1\over |\xi|}\right)^{\left({n+1\over n}\right)\alpha}\cos\left(2\pi|\xi|+{\pi\over 2}\left({n+1\over n}\right)\alpha\right)d\xi\right\}du\right|^q dx \right\}^{1\over q}  
\\ \ds~~~~~~~~~~~~~~~~~~~~~~~~~~~~~~~~~~~~~~~~~~~~~~~~~~~~~~~~~~~~~~~~~~~~~~~~~~~~~~~~~~~~~~~~~~~~~~~~~~~~~~~~~~~~~~ u\mt (r-s) u
\\ \ds
~\leq~\C_\alpha\left({1\over rs}\right)^{\left({n+1\over 2n}\right)\alpha}(r-s)^{\left({n+1\over n}\right)\alpha+{n\over q}}\left\{\int_{\R^n} \left|f\left((r-s)x\right)\right|^{q\over q-1} dx \right\}^{q-1\over q}  \qquad\hbox{\small{by assumption}}
\\\\ \ds
 ~=~\C_\alpha \left({1\over rs}\right)^{\left({n+1\over 2n}\right)\alpha}\left({1\over r-s}\right)^{n\left({q-1\over q}\right)-{n\over q}-\left({n+1\over n}\right)\alpha}  
 \left\{\int_{\R^n} \left|f(x)\right|^{q\over q-1} dx \right\}^{q-1\over q}  \qquad x\mt (r-s)^{-1}x
 \\\\ \ds
 ~=~\C_\alpha \left({1\over rs}\right)^{\left({n+1\over 2n}\right)\alpha}\left({1\over r-s}\right)^{\left({n-1\over n}\right)\alpha}   \left\{\int_{\R^n} \left|f(x)\right|^{q\over q-1} dx \right\}^{q-1\over q}. 
  \end{array}
 \eeq
 A vice versa estimate works for $r-s$ replaced by $r+s$.

Let $\Omega^{\lambda(\alpha)}$  be defined by the inverse Fourier transform of (\ref{Omega^alpha Transform}) for $0<\alpha<n$. Recall (\ref{convolution inequality})  from {\bf Proposition 2}. We have
  \bel{Est z} 
 \left\| f\ast\Omega^{\lambda(2\alpha)}\right\|_{\L^q(\R^n)}~\leq~\C_{\alpha~v}~\left\| f\right\|_{\L^{q\over q-1}(\R^n)}
 \eeq  
for
 \bel{Range 2alpha}
 {q-1\over q}~=~{1\over 2}\left(1+{2\alpha\over n}\right),\qquad {1\over q}~=~{1\over 2}\left(1-{2\alpha\over n}\right). 
\eeq
From (\ref{Omega split})-(\ref{Omega split S+E}), we have $\Hat{\Omega}^{\lambda(2\alpha)}=\Hat{\mathfrak{S}}^{\lambda(2\alpha)}(\xi)+\Hat{\mathcal{E}}^{\lambda(2\alpha)}$ where
 \bel{E 2alpha Est}
 \begin{array}{lr}\ds
 \left|\Hat{\mathcal{E}}^{\lambda(2\alpha)}\right|~=~\left|  \left({1\over |\xi|}\right)^{\left({n+1\over n}\right)\alpha-{1\over 2}+\i v}\e_v\left(2\pi|\xi|\right) \right| 
 \\\\ \ds~~~~~~~~~~~~
~\leq~   \left({1\over |\xi|}\right)^{\left({n+1\over n}\right)\alpha-{1\over 2}}\left|\e_v\left(2\pi|\xi|\right) \right|
~\leq~\C~e^{\pi|v|}~\left({1\over |\xi|}\right)^{2\alpha}
 \end{array}
 \eeq
can be deduced from  (\ref{E norm est 1})-(\ref{E norm est 3}). 

Hardy-Littlewood-Sobolev theorem \cite{Hardy-Littlewood}-\cite{Sobolev} implies that $\Hat{\mathcal{E}}^{\lambda(2\alpha)}(\xi)$ is an $\L^{q\over q-1}(\R^n)\mt\L^q(\R^n)$ Fourier multiplier. Together with  (\ref{Est z})-(\ref{Range 2alpha}), we must have 
\bel{S 2alpha}
\Hat{\mathfrak{S}}^{\lambda(2\alpha)}(\xi)~=~{1\over \pi}  \left({1\over |\xi|}\right)^{\left({n+1\over n}\right)\alpha+\i v} \cos\left(2\pi|\xi|-{\pi\over 2}\left({n+1\over n}\right)\alpha-{\pi\over 2}(1+\i v)\right)
\eeq
to be an $\L^{q\over q-1}(\R^n)\mt\L^q(\R^n)$ Fourier multiplier. In particular, at $v=0$, (\ref{S 2alpha})  shows that 
$\left({1\over |\xi|}\right)^{\left({n+1\over n}\right)\alpha}\sin\left(2\pi|\xi|+{\pi\over 2}\left({n+1\over n}\right)\alpha\right)$
is an $\L^{q\over q-1}(\R^n)\mt\L^q(\R^n)$  Fourier multiplier as desired. 

Note that $|\xi|^{\i v}$ is a $\L^p$ Fourier multiplier for $1<p<\infty$ at any $v\in\R$. See p. 51 of Stein \cite{Stein'}. From (\ref{S 2alpha}), we also have
\bel{S* 2alpha}
\mathfrak{S}^\ast_v(\xi)~\doteq~ \left({1\over |\xi|}\right)^{\left({n+1\over n}\right)\alpha} \cos\left(2\pi|\xi|-{\pi\over 2}\left({n+1\over n}\right)\alpha-{\pi\over 2}(1+\i v)\right)
\eeq
as an $\L^{q\over q-1}(\R^n)\mt\L^q(\R^n)$ Fourier multiplier. 
Moreover, for every $v\in\R$,  the identity 
\bel{identity}
\begin{array}{lr}\ds
~~~~~~~\cos\left(2\pi|\xi|-{\pi\over 2}\left({n+1\over n}\right)\alpha-{\pi\over 2}(1+\i v)\right)
\\\\ \ds
~=~\sin\left(2\pi|\xi|-{\pi\over 2}\left({n+1\over n}\right)\alpha\right)\cos\left({\i\pi v\over 2}\right)-\cos\left(2\pi|\xi|-{\pi\over 2}\left({n+1\over n}\right)\alpha\right)\sin \left({\i\pi v\over 2}\right)
\end{array}
\eeq
further implies that 
$\left({1\over |\xi|}\right)^{\left({n+1\over n}\right)\alpha}\cos\left(2\pi|\xi|+{\pi\over 2}\left({n+1\over n}\right)\alpha\right)$ is another desired  Fourier multiplier.

\section{Proof of Theorem One}
\setcounter{equation}{0}
Let  $0<\alpha\leq{n\over n+1}$.
Recall {\bf Proposition 1} and {\bf 2} from section 2. The estimate (\ref{Result one}) implies 
\bel{mix-norm.1}
\begin{array}{cc}\ds
  \left\| |r|^\alpha \left\| f(\cdot,t)\ast\Omega^{\alpha}_{|r|}\right\|_{\L^{q_1}(\R^n)}\right\|_{\L^{q_1}(\R,~ |r|^{-1}dr)}~\leq~\C_{\alpha} \left\| f(\cdot,t)\right\|_{\L^2(\R^n)},\qquad t\in\R
  \\\\ \ds
  \hbox{for}\qquad  {\alpha\over n}~=~{1\over 2}-{1\over q_1},\qquad 1<q_1<\infty.  
 \end{array}
 \eeq  
By changing dilations in (\ref{convolution inequality}) at $v=0$, we have
 \bel{Est1}
 \begin{array}{cc}\ds
  \left\| |r|^\alpha\Big(f(\cdot,t)\ast\Omega^{\alpha}_{|r|}\Big)\right\|_{\L^{q_o}(\R^n)}~\leq~\C_{\alpha}~\left\| f(\cdot,t)\right\|_{\L^{p_o}(\R^n)}\qquad t\in\R
\\\\ \ds
\hbox{for}\qquad {1\over p_o}~=~{1\over 2}\left(1+{\alpha\over n}\right),\qquad {1\over q_o}~=~{1\over 2}\left(1-{\alpha\over n}\right).
\end{array}
\eeq
By taking the supremum over $r\neq0$, (\ref{Est1}) can be interpreted as
\bel{mix-norm.2}
 \begin{array}{cc}\ds
\left\| |r|^\alpha \left\| f(\cdot,t)\ast\Omega^{\alpha}_{|r|}\right\|_{\L^{q_o}(\R^n)}\right\|_{\L^\infty(\R,~ |r|^{-1}dr)}~\leq~\C_{\alpha} \left\| f(\cdot,t)\right\|_{\L^{p_o}(\R^n)},\qquad t\in\R
\\\\ \ds
\hbox{for}\qquad {1\over p_o}~=~{1\over 2}\left(1+{\alpha\over n}\right),\qquad {1\over q_o}~=~{1\over 2}\left(1-{\alpha\over n}\right).
\end{array}
\eeq 
Let $0<\theta<1$ and 
 \bel{s}
 \begin{array}{cc}\ds
 {1\over s}~\doteq~{\theta\over q_1},
\qquad
 {1\over p}~=~{1-\theta\over p_o}+{\theta\over 2},\qquad {1\over q}~=~{1-\theta\over q_o}+{\theta\over q_1}.
 \end{array}
 \eeq
 Note that 
 \bel{constraint+condition}
 \begin{array}{ccc}\ds
 s~>~q, \qquad {1\over 2}~<~{1\over p}~<~{1\over 2}+{\alpha\over 2n} 
  \end{array}
 \eeq
 because of (\ref{Est1}) and $0<\theta<1$. 

 From (\ref{mix-norm.1})-(\ref{constraint+condition}), by applying Riesz-Thorin interpolation theorem in mixed-norms, we have
 \bel{Est2}
 \begin{array}{cc}\ds
 \left\{\int_\R |r|^{\alpha s}\left\{\int_{\R^n}  \left|\Big( f(\cdot,t)\ast\Omega^{\alpha}_{|r|}\Big)(x)\right|^qdx\right\}^{s\over q} {dr\over |r|}\right\}^{1\over s}~\leq~\C_{\alpha} \left\| f(\cdot,t)\right\|_{\L^p(\R^n)},\qquad t\in\R
 \\\\ \ds
 \hbox{if}\qquad {\alpha\over n}~=~{1\over p}-{1\over q},\qquad {1\over 2}~<~{1\over p}~<~{1\over 2}+{\alpha\over 2n}.
 \end{array}
 \eeq

 Recall $\I_{\alpha}$ defined in (\ref{I_alpha}). 
 Given $f\in\L^p(\R^{n+1})$ and $g\in\L^{q\over q-1}(\R^{n+1})$,  we have
 \bel{Bilinear form}
\begin{array}{lr}\ds
\iint_{\R^{n+1}}\Big(\I_{\alpha}f\Big)(x,t)  g(x,t)dxdt~=~
\iint_{\R^{n+1}}
\left\{\int_\R \Big(f(\cdot,t-r)\ast\Omega^{\alpha}_{|r|}\Big)(x)| r|^{\left({n+1\over n}\right)\alpha-1}dr\right\}g(x,t)dxdt
\\\\ \ds
~\leq~\iint_{\R^{2}}\left|\int_{\R^n} \Big(f(\cdot,t)\ast\Omega^{\alpha}_{|r|}\Big)(x)g(x,t+r)dx\right|  | r|^{\left({n+1\over n}\right)\alpha-1}drdt  \qquad t\mt t+r
\\\\ \ds
~\leq~\int_\R\int_{\R} \left\|f(\cdot,t)\ast\Omega^{\alpha}_{|r|} \right\|_{\L^q(\R^{n})}\left\|g(\cdot, t+r)\right\|_{\L^{q\over q-1}(\R^{n})} |r|^{\left({n+1\over n}\right)\alpha-1}dr dt
\\\\ \ds
~\leq~\int_\R \left\{\int_{\R} \left\|f(\cdot,t)\ast\Omega^{\alpha}_{|r|}\right\|^s_{\L^q(\R^{n})} |r|^{\alpha s} {dr\over |r|}\right\}^{1\over s}
\left\{\int_\R  \left\|g(\cdot, t+r)\right\|_{\L^{q\over q-1}(\R^{n})}^{s\over s-1} |r|^{\left(\left({n+1\over n}\right)\alpha-1+{1\over s}-\alpha\right)\left({s\over s-1}\right)} dr\right\}^{s-1\over s} dt
\\\\ \ds
~\leq~\left\{\int_\R \left\{\int_{\R} \left\|f(\cdot,t)\ast\Omega^{\alpha}_{|r|} \right\|^s_{\L^q(\R^{n})} |r|^{\alpha s} {dr\over |r|}\right\}^{p\over s} dt\right\}^{1\over p}
\\\\ \ds~~~~~~~
\left\{\int_\R\left\{\int_\R  \left\|g(\cdot, t+r)\right\|_{\L^{q\over q-1}(\R^{n})}^{s\over s-1} |r|^{\left(\left({n+1\over n}\right)\alpha-1+{1\over s}-\alpha\right)\left({s\over s-1}\right)} dr\right\}^{\left({s-1\over s}\right)\left({p\over p-1}\right)} dt\right\}^{p-1\over p}
\end{array}
\eeq
by using H\"{o}lder inequality 3 times.

Direct computation shows
\bel{Formula computa}
\begin{array}{cc}\ds
\left(\left(n+1\right)\left({\alpha\over n}\right)-1+{1\over s}-\alpha\right)\left({s\over s-1}\right)~=~\left({s\over s-1}\right)\left({\alpha\over n}\right)-1
\\\\ \ds
\hbox{where}\qquad {s\over s-1}~<~{q\over q-1},\qquad {\alpha\over n}~=~{q-1\over q}- {p-1\over p}.
\end{array}
\eeq 
From (\ref{Formula computa}), by applying Hardy-Littlewood-Sobolev inequality \cite{Hardy-Littlewood}-\cite{Sobolev} on $\R$, we have
\bel{1-d HLS est}
\begin{array}{lr}\ds
\left\{\int_\R\left\{\int_\R  \left\|g(\cdot, t+r)\right\|_{\L^{q\over q-1}(\R^{n})}^{s\over s-1} |r|^{\left(\left(n+1\right)\left({\alpha\over n}\right)-1+{1\over s}-\alpha\right)\left({s\over s-1}\right)} dr\right\}^{\left({s-1\over s}\right)\left({p\over p-1}\right)} dt\right\}^{p-1\over p}
\\\\ \ds
~\leq~\C_{p~q}~\left\{\int_\R \left\|g(\cdot, t)\right\|_{\L^{q\over q-1}(\R^{n})}^{q\over q-1}dt\right\}^{q-1\over q}~=~\C_{p~q}~\left\|g\right\|_{\L^{q\over q-1}(\R^{n+1})}.
\end{array}
\eeq
On the other hand,  (\ref{Est2}) implies
\bel{critical Est} 
\begin{array}{lr}\ds
\left\{\int_\R \left\{\int_{\R} \left\|f(\cdot,t)\ast\Omega^{\alpha}_{|r|} \right\|^s_{\L^q(\R^{n})} |r|^{\alpha s} {dr\over |r|}\right\}^{p\over s} dt\right\}^{1\over p}~\leq~\C_{p~q}~\left\| f\right\|_{\L^p(\R^{n+1})}.
\end{array}
\eeq
By bringing (\ref{1-d HLS est}) and (\ref{critical Est}) back to (\ref{Bilinear form}) and taking the supremum over $\left\|g\right\|_{\L^{q\over q-1}(\R^{n+1})}=1$, we have
\bel{I_lambda regularity}
\begin{array}{cc}\ds
\left\|\I_\alpha f\right\|_{\L^q(\R^{n+1})}~\leq~\C_{p~q}~\left\| f\right\|_{\L^p(\R^{n+1})}, \qquad 1<p<q<\infty
\\\\ \ds
\hbox{if}\qquad {\alpha\over n}~=~{1\over p}-{1\over q},\qquad {1\over 2}~<~{1\over p}~<~{1\over 2}+{\alpha\over 2n}.
\end{array}
\eeq
Lastly, it is easy to verify that  the adjoint operator of $\I_{\alpha}$ satisfies all the regarding estimates. By taking into account  
\bel{range computa.1-2}
\begin{array}{cc}\ds
{1\over 2}-{\alpha\over n}~<~{1\over q}~<~{1\over 2}-{\alpha\over 2n}\qquad\Longleftrightarrow\qquad {1\over 2}+{\alpha\over 2n}~<~{q-1\over q}~<~{1\over 2}+{\alpha\over n},
\\\\ \ds
{1\over 2}~<~{1\over p}~<~{1\over 2}+{\alpha\over 2n}\qquad\Longleftrightarrow\qquad {1\over 2}-{\alpha\over 2n}~<~{p-1\over p}~<~{1\over 2}
\end{array}
\eeq
and applying Riesz-Thorin interpolation theorem again, we obtain (\ref{RESULT*}) as desired.

\v
{\small Department of Mathematics, Westlake University}

{\small email: wangzipeng@westlake.edu.cn}

\end{document}